# The hardest logic puzzle ever becomes even tougher

Abstract

"The hardest logic puzzle ever" presented by George Boolos became a target for philosophers and logicians who tried to modify it and make it even tougher. I propose further modification of the original puzzle where part of the available information is eliminated but the solution is still possible. The solution also gives interesting ideas on logic behind discovery of unknown language.

## 1. Introduction

"The hardest logic puzzle ever," originally presented by Boolos (1996), had already been amended several times in order to make it tougher. Rabern and Rabern (2008) have modified the behavior of Random to make it really random and avoid trivialization. Rabern and Rabern (2008) and Uzquiano (2010) analyzed a two-question solution using so called "head-exploding questions".

The purpose of this note is to offer further modification that clearly makes the puzzle tougher. The solution for this modified puzzle explores our ability to extract information in situations of complete language ignorance by exploiting features of any possible language.

## 2. Original puzzle and modifications

Let's first recall Boolos' statement of the puzzle:

*"Three gods A, B, and C are called, in some order, True, False, and Random. True always speaks truly, False always speaks falsely, but whether Random speaks truly or falsely is a completely random matter. Your task is to determine the identities of A, B, and C by asking three yes-no questions; each question must be put to exactly one god. The gods understand English, but will answer all questions in their own language in which the words for 'yes' and 'no' are 'da' and 'ja', in some order. You do not know which word means which."*

*Clarifications:*

*(1) It could be that some god gets asked more than one question (and hence that some god is not asked any question at all).*

*(2) What the second question is, and to which god it is put, may depend on the answer to the first question. (And of course, similarly for the third question)*

*(3) Whether Random speaks truly or not should be thought of as depending on the flip of a coin hidden in his brain: if the coin comes down heads, he speaks truly; if tails, falsely.*

*(4) Random will answer 'da' or 'ja' when asked any yes-no question.*

If Random speaks the truth or lies randomly, there is a very simple solution for this puzzle. Rabern and Rabern (2008) have proposed to modify the third point above with the following:

*"Whether Random answers 'da' or 'ja' should be thought of as depending on the flip of a coin hidden in his brain: if the coin comes down heads, he answers 'yes'; if tails, 'no'."*

Rabern and Rabern (2008) and Uzquiano (2010) showed that one can solve the puzzle in two questions using self-referential questions. Uzquiano also gave two-question solution without self-reference, by restricting the knowledge of the gods concerning Random's behavior. Now, we can either require a two-question solution for the puzzle or add a clarification:

*"If you ask a question, where the god can't give an answer, he will answer 'no'."*

### 3. New modification: Complete language ignorance

It sounds very unnatural that we know the words 'da' and 'ja' and don't know their meaning (if a situation where we speak with three gods can sound natural). Let's assume that we don't know anything about their language. We know that gods speak an unknown language and will answer 'yes' or 'no' in their language. We should also assume that all gods speak the same language and will use the same words (not synonyms). This modification makes the solution much tougher.

Now there are three independent modifications of the puzzle that can be combined: random behavior of Random, paradoxical questions permitted and complete language ignorance. This article analyzes three statements of the puzzle in case of complete language ignorance:

I. Original behavior of Random, paradoxical questions prohibited, provided solution in three questions

II. Random behaves really randomly, paradoxical questions prohibited, provided proof of no solution in three questions

III. Random behaves really randomly, paradoxical questions permitted, provided solution in two questions

Why are these modifications tougher? Normally this is difficult to compare two puzzles by level of "hardness". But it is possible in case, if two puzzles are similar and all solutions for one are a subset of all solutions for another. For example, puzzle A has two solutions 'a' and 'b'. Puzzle B has only one solution 'a'. Than puzzle B is clearly tougher, because it is harder (or at least not easier) to find the solution. For the sake of rigorousness, we should also exclude the possibility that formulation of one puzzle provides hints for the solution.

Puzzle II above is tougher than any 3 question versions with paradoxical questions prohibited, because:

1. Puzzle II above has only one solution that is provided later in this article.
2. Original puzzle has 3 solutions:
   (a) Original solution of Boolos
   (b) 2 question solution proposed by Rabern and Rabern (2008) that exploits not fully random behavior of Random
   (c) Solution for puzzle II that also works with original puzzle.

3. Modification by Raberns (Random behaves really randomly) has 2 solutions:

    (a) Original solution of Boolos

    (b) Solution for puzzle II that also works with this puzzle.

Puzzle III above is tougher than any 2 question formulations with paradoxical questions permitted, because:

1. Puzzle III has only one solution that is provided later in this article.

2. Modification by by Rabern and Rabern (2008) has 3 solutions:

    (a) Solution using self-referential questions

    (b) Solution by restricting the knowledge of the gods concerning Random's behavior, Uzquiano (2010)

    (c) Solution for puzzle III that also works with this puzzle.

4. *Solution for Puzzle I: Original behavior of Random, paradoxical questions prohibited, solution in three questions*

Full formulation of amended puzzle:

*"Three gods A, B, and C are called, in some order, True, False and Random. True always speaks truly, False always speaks falsely, but whether Random speaks truly or falsely is a completely random matter. Your task is to determine the identities of A, B, and C by asking three yes/no questions; each question must be put to exactly one god. The gods understand English, but will answer all questions in their own language. You do not know their language."*

*Clarifications:*

   *(1) "It could be that some god gets asked more than one question (and hence that some god is not asked any question at all).*

*(2) What the second question is, and to which god it is put, may depend on the answer to the first question. (And of course, similarly for the third question)*

*(3) Whether Random speaks truly or not should be thought of as depending on the flip of a coin hidden in his brain: if the coin comes down heads, he speaks truly; if tails, falsely.*

*(4) If you ask a question, where the god can't give an answer, he will answer 'no'.*

*(5) All gods speak in the same language and will use the same words (not synonyms)."*

Before continuing with this article, the reader may wish to pause and attempt a solution.

Since there is no information about the gods' answers, it is impossible to refer to them in questions. Knowing only one word ("yes" or "no") would be enough. Nevertheless, we have some information we can use for reference. Before asking questions, we have only one piece of information—gods have a language. It stands to reason that they must have two clearly different words for "yes" or "no." We can exploit this. First, we ask the god to sort these two words by a sorting rule, and second, we refer to the "first" word using an *embedded question*, introduced by Rabern and Rabern (2008).

A simple sorting rule could be *alphabetical order in English transliteration.* However, we don't know anything about the god's language. Maybe they will answer by signs or produce different tones, etc. Subsequently, we need a *universal sorting rule*.

*Universal Sorting Rule lemma.* There is a sorting rule that can allow us to determine the order of "yes" or "no" in any possible language.

*Proof.* Consider the following rule:

"Sort in alphabetical order, descriptions that I will give to words 'yes' and 'no' in your language once you tell them to me."

1. Gods can foresee what description the person who asks will give to the words "yes" and "no," in their language[1].

2. Hence, there is a difference in these words or signs that could be perceived by the person who asks; this person will produce different descriptions.

3. These descriptions will be given in English and could be sorted in alphabetical order.

The universal sorting rule allows us to solve a given puzzle with the following questions. The first two questions are predetermined:

(1st question, to A): "If I asked you, in your current mental state, 'Are you Random?,' would you answer with the word which comes first alphabetically in the list of descriptions that I would give to the words 'yes' and 'no' in your language, assuming you had told them to me previously."

(2nd question, to B): "If I asked you, in your current mental state, 'are you Random?,' would you answer with the same response given to my previous question?"

After these two questions, we can get three outcomes:

(1) Some answer (meaning inconclusive)—Same answer (Yes): We know that B answered Yes and he is Random.

(2) Some answer (Yes)—Different answer, second when sorted (No): We know that A responded Yes and he is Random.

---

[1] In this proof, we rely on the god's ability to predict our behavior. Alternatively, we can communicate a sound and vision recognition algorithm in our question that could be used to generate descriptions. If we believe in progress in Artificial Intelligence development, it is a possible option in principle.

(3) Some answer (No)—Different answer, first when sorted (No): We know that both gods answered No and C is Random.

Now we know who is Random and can define True or False by asking one of them, "If I asked you, 'are you True?,' would you answer with the response given to my previous question?" The same answer as previously means that you asked True, otherwise you asked False.

5. *Solution for Puzzle II: Real Random, prohibited paradoxical questions—in three questions*

Let's consider a modification of the puzzle that includes complete language ignorance; Random answers depend only on the coin flip and head-exploding questions are prohibited.

It seems that in this situation, there is no solution in three questions. At least we can prove that it is impossible using only the sorting rule. If there is a solution, one should invent a way to extract more information from god's answers than with the sorting rule.

*Proof*

1. Using the sorting rule we can get seven different outcomes after three questions:

    (1) Some answer—Same answer—Same answer

    (2) Some answer—Same answer—Answer first when sorted

    (3) Some answer—Same answer—Answer second when sorted

    (4) Some answer—Answer second when sorted—Answer first when sorted

    (5) Some answer—Answer second when sorted—Answer second when sorted

    (6) Some answer—Answer first when sorted—Answer first when sorted

    (7) Some answer—Answer first when sorted—Answer second when sorted

The only way we can solve the puzzle is to allocate every possible order of gods to only one outcome using our questions. There are only six possible ways to order the gods, so the puzzle seems solvable. However, a problem arises when we consider random answers of Random.

2. We can't extract any information from the first question alone, so the first two questions are predetermined in any possible strategy.

3. A successful strategy can't address the first two questions to one god as this god may be Random. It follows that after two questions, we will have three possible outcomes, and in each case, we can't exclude the possibility that A is Random:

   (1) Some answer—Same answer (3 outcomes to third question): RTF[2], RFT

   (2) Some answer—Another answer second when sorted (two outcomes to the third question): RTF, RFT

   (3) Some answer—Another answer first when sorted (two outcomes to the third question): RTF, RFT

Using different questions, we can allocate TFR, FTR, TRF, FRT to some of the outcomes. But there is only one place left, so there is no solution in this case.

4. If we address the first two questions to A and B, there is a possibility that one of them is Random. After two questions, we can encounter 10 different possibilities and this doesn't depend on what exact questions we ask:

   (1) TFR

   (2) FTR

   (3) R(answers yes)TF

   (4) R(answers no)TF

---

[2] Here and later RTF means Random, True, False

(5) R(answers yes)FT

(6) R(answers no)FT

(7) TR(answers yes)F

(8) TR(answers no)F

(9) FR(answers yes)T

(10) FR(answers no)T

5. Each of the above cases should be allocated to possible outcomes. Only the outcome *Some answer—Same answer—Same answer* can include two cases with the same order of gods (like R(1)TF and R(2)TF), hence we don't need to differ between them. However, all cases with the same order of gods can't be allocated by any other single outcome, because they have different answers to questions, by definition.

6. As a result, we have only seven outcomes for nine cases. No strategy can provide us with a solution.

6. *Solution for Puzzle III: Random behaves really randomly, paradoxical questions permitted, solution in two questions*

Rabern and Rabern (2008) use the image of head explosion to show the case when the god can't answer Yes or No to the question. For our case we should assume that if we see a head explosion, we can clearly say that it is a head explosion, not a sign of "yes" or "no." It should also be clarified that Random's head can't explode because he answers randomly to any question. There might be also confusion whether you can ask the same god after his head "exploded". My assumption is - you can't. The solution is still possible, but conditions in this case are tougher.

Here is the solution. First question directed to A:

(1*) Would you answer "with the word which comes first alphabetically in the list of descriptions..." to the question whether either:

(a) B isn't Random and you are False, or

(b) B is Random and you would answer "with the word which comes second alphabetically in the list of descriptions…" to (1*)?

Uzquiano (2010) showed that this question has the following outcomes:

(1) If B is Random—Explosion

(2) If A is Random—Random answer

(3) If C is Random and A is False—first answer when sorted

(4) If C is Random and A is True—second answer when sorted

Second question:

After the first question, we have two possible outcomes—head explosion or some answer. It will influence our strategy.

If the first answer was explosion, then ask C:

(2*) Would you answer "with the word which comes first alphabetically in the list of descriptions…" to the question whether true is both statements:

(a) You are True, and

(b) You would answer "with the word which comes second alphabetically in the list of descriptions…" to (2*)?

If the third god is True, his answer is Explosion. We now know the identities of the gods to be False, Random, and True.

If the third god is False, his answer is second when sorted. Actually, we cannot define whether the answer is first or second when sorted, but we already know that gods' order is True, Random, and False.

If the first answer was not explosion, then ask B:

(2*) Would you answer "with the response given to my previous question" to the question whether either:

(a) A is Random, you are True, and you would answer "with the response opposite to the one given to my previous question" to (2*), or

(b) A is Random, and you are False?

If A isn't Random, the embedded statement is False and B will answer unlike A. We can now identify which answer is first and second when sorted. If A is False, then his answer should be first when sorted and the right order of gods is False, True, and Random. If A is True, then his answer should be second when sorted and the right order of gods is True, False, and Random.

If A is Random and B is True, then B cannot answer and his head explodes. We can conclude that the right order of gods is Random, True, and False.

If A is Random and B is False, then the answer will be the same as given by A. We can conclude that the right order of gods is Random, False, and True.

*Nikolay Novozhilov, nikolay.novozhilov@gmail.com*